\documentstyle[12pt,draft]{article}

 \newcounter{abceqn}

 \newcounter{abcfig}

\renewcommand{\H}{{\cal H}}
\newcommand{\E}{{\cal E}}

\newcommand{\e}{\varepsilon}

\newcommand{\tp}{\tilde{p}}
\newcommand{\tq}{\tilde{q}}

\newcommand{\k}{\kappa}
\newcommand{\ga}{\gamma}

\newcommand{\dl}{\delta}
\newcommand{\Dl}{\Delta}
\newcommand{\th}{\theta}
\newcommand{\vth}{\vartheta}
\newcommand{\vphi}{\varphi}

\newcommand{\ra}{\rightarrow}
\newcommand{\al}{\alpha}
\newcommand{\be}{\beta}
\newcommand{\sg}{\sigma}
\newcommand{\Sg}{\Sigma}

\newcommand{\pa}{\partial}

\newcommand{\hq}{\hat{q}}
\newcommand{\hp}{\hat{p}}

\newcommand{\bh}{\bar{h}}

\newcommand{\la}{\lambda}
\newcommand{\bq}{\bar{q}}
\newcommand{\bp}{\bar{p}}

\newcommand{\nid}{\noindent}

\newcommand{\om}{\omega}
\newcommand{\Om}{\Omega}

\renewcommand{\theequation}{\thesection.\arabic{equation}}
\newcommand{\eqnsection}[1]{
	\section{#1}
	\setcounter{equation}{0}
	\renewcommand{\theequation}{\thesection.\arabic{equation}}
	\setcounter{figure}{0}
	\renewcommand{\thefigure}{\thesection.\arabic{figure}}
	\setcounter{remark}{0}
	\renewcommand{\theremark}{\thesection.\arabic{remark}}
	\setcounter{theorem}{0}
	\renewcommand{\thetheorem}{\thesection.\arabic{theorem}}
	\setcounter{lemma}{0}
	\renewcommand{\thelemma}{\thesection.\arabic{lemma}}
	\setcounter{proposition}{0}
	\renewcommand{\theproposition}{\thesection.\arabic{proposition}}
}

\title{\bf Singularly Perturbed Vector and Scalar Nonlinear Schr\"odinger 
Equations with Persistent Homoclinic Orbits}

\author{\\ \\ \\ \\ Yanguang\ (Charles) \ \ Li  
\thanks{This work is partially supported by a Guggenheim Fellowship.}
\\ \\ \\ Department of Mathematics \\ \\ University of Missouri - Columbia
\\ \\ Columbia, MO 65211 \\ \\ E-mail: cli@math.missouri.edu}

\date{\today}

\renewcommand{\theequation}{\thesection.\arabic{equation}}

\begin{document}
\bibliographystyle{plain}
\maketitle

\newpage
\begin{abstract} 
Singularly perturbed vector nonlinear Schr\"odinger 
equations (PVNLS) are investigated. Emphasis is placed 
upon the relation with their restriction: The singularly 
perturbed scalar nonlinear Schr\"odinger equation (PNLS)
studied in \cite{Li01a}. It turns out that the persistent 
homoclinic orbit for the PNLS \cite{Li01a} is the only one 
for the PVNLS, asymptotic to the same saddle.
\end{abstract}

\newtheorem{lemma}{Lemma}
\newtheorem{theorem}{Theorem}
\newtheorem{corollary}{Corollary}
\newtheorem{remark}{Remark}
\newtheorem{definition}{Definition}
\newtheorem{proposition}{Proposition}
\newtheorem{assumption}{Assumption}

\newpage
\tableofcontents

\newpage
\eqnsection{Introduction}

In recent years, novel results have been obtained on the solutions of 
the vector nonlinear Schr\"odinger equations \cite{AOT99} \cite{AOT00} 
\cite{YT01}. Abundant ordinary integrable results have been carried through 
\cite{WF00} \cite{FSW00}, including linear stability calculations 
\cite{FMMW00}. Specifically, the vector nonlinear 
Schr\"odinger equations can be written as
\begin{eqnarray*}
& & ip_t + p_{xx} + \frac{1}{2} (|p|^2 + \chi |q|^2) p = 0 , \\
& & iq_t + q_{xx} + \frac{1}{2} (\chi |p|^2 + |q|^2) q = 0 , 
\end{eqnarray*}
where $p$ and $q$ are complex valued functions of the two real variables 
$t$ and $x$, and $\chi$ is a positive constant. These equations describe
the evolution of two orthogonal pulse envelopes in birefringent optical 
fibers \cite{Men87} \cite{Men89}, with industrial applications in fiber 
communication systems \cite{HK95} and all-optical switching devices 
\cite{Isl92}. For linearly birefringent fibers \cite{Men87}, $\chi =2/3$.
For elliptically birefringent fibers, $\chi$ can take other positive 
values \cite{Men89}. When $\chi = 1$, these equations are first 
shown to be integrable by S. Manakov \cite{Man74}, and thus called Manakov 
equations. When $\chi$ is not 1 or 0, these equations are 
non-integrable. Propelled by the industrial applications, extensive 
mathematical studies on the vector nonlinear Schr\"odinger equations 
have been carried. This article is another piece of mathematical works.
Like the scalar nonlinear Schr\"odinger equation, the vector nonlinear 
Schr\"odinger equations also possess figure eight structures in their 
phase space. Figure eight structures are also called separatrices. 
For two-dimensional Hamiltonian systems, figure eight structures 
are given by the singular level sets of the Hamiltonian.
Our goal is to understand the consequence of such figure 
eight structures when these equations are under perturbations, in 
particular, whether or not the consequence is chaotic dynamics. 
Specifically, we will investigate the type of chaos created through 
homoclinic orbits. (cf: the section ``Conclusion and Discussion''.)

This article will deal with the problem on the existence of homoclinic 
orbits for the singularly perturbed vector nonlinear Schr\"odinger equations,
\begin{eqnarray}
& & ip_t + p_{xx} + \frac{1}{2}(|p|^2 + |q|^2) p = 
i \e [ p_{xx} -\al p - \be e^{i \frac{1}{2}\om^2 t}]\ , \label{inls1}\\
& & iq_t + q_{xx} + \frac{1}{2} (|p|^2 + |q|^2) q = 
i \e [ q_{xx} -\al q - \be e^{i \frac{1}{2}\om^2 t}]\  , \label{inls2}
\end{eqnarray}
where $p(t,x)$ and $q(t,x)$ are subject to periodic boundary condition 
of period $2\pi$, and are even in $x$, i.e. 
\[
p(t,x + 2 \pi) = p(t,x)\ , \ \ p(t,-x) = p(t,x)\ , 
\]
\[
q(t,x + 2 \pi) = q(t,x)\ , \ \ q(t,-x) = q(t,x)\ , 
\]
$\om \in (1,2)$, $\al > 0$ and $\be$ are real constants, and 
$\e > 0$ is the perturbation parameter. 

By the transformation,
\[
p = \tilde{p} e^{i \frac{1}{2} \om^2 t}\ , \ \ 
q = \tilde{q} e^{i \frac{1}{2} \om^2 t}\ , 
\]
and dropping $\ \tilde{}\ $, one gets an autonomous system,
\begin{eqnarray}
& & ip_t + p_{xx} + \frac{1}{2} [(|p|^2 + |q|^2)-\om^2] p = 
i \e [ p_{xx} -\al p - \be ]\ , \label{pnls1}\\
& & iq_t + q_{xx} + \frac{1}{2} [(|p|^2 + |q|^2)-\om^2] q = 
i \e [ q_{xx} -\al q - \be ]\  , \label{pnls2}
\end{eqnarray}
The rest of this article will be dealing with this form of 
perturbed vector nonlinear Schr\"odinger equations (PVNLS).

The corresponding work 
for regularly perturbed scalar nonlinear Schr\"odinger equation 
has been accomplished in \cite{LMSW96}, where the singular perturbation 
$\e \pa_x^2$ is replaced by a regular perturbation $\e \hat{\pa}_x^2$ 
which is a bounded Fourier multiplier obtained from a truncation 
of $\e \pa_x^2$. The main difficulty introduced by the singular 
perturbation $\e \pa_x^2$ is that it breaks the spectral gap condition 
of the unperturbed system. Thus the standard invariant manifold results 
do not hold. Nevertheless, certain invariant manifold results
can be established, which are sufficient for establishing the existence 
of homoclinic orbits. The corresponding work for such singularly 
perturbed scalar nonlinear Schr\"odinger equation (PNLS) has been 
accomplished in 
\cite{Li01a}. The crucial new features of this article are listed below:
\begin{enumerate}
\item Homoclinic orbits asymptotic to the spatially uniform and 
temporally periodic 
solutions up to phase shifts, are constructed for the integrable vector 
nonlinear Schr\"odinger equations ((\ref{inls1}, \ref{inls2}) at $\e =0$). 
In particular, such homoclinic orbits are built with spatially non-periodic 
and non-anti-periodic Bloch eigenfunctions. This study provides an elegant 
formula which is missing in \cite{WF00} \cite{FSW00} \cite{FMMW00}.
\item For the scalar nonlinear Schr\"odinger equation, a 
complete sequence of Melnikov vectors can be constructed out of the 
Floquet discriminant for the Lax pair \cite{LM94}. 
Unlike scalar nonlinear Schr\"odinger equation, 
the Lax pair for the integrable vector nonlinear Schr\"odinger 
equations, does not have a complete isospectral theory built upon
a Floquet discriminant. 
By restricting $\om$ in the interval ($1,2$), we are able to deal 
with one unstable mode problem for which only one Melnikov function is needed. The Melnikov function will be built 
out of the Hamiltonian and the invariant $L^2$ norms. 
\item In the 4D invariant subspace of functions independent of $x$, 
the dynamics is quite different from that of perturbed scalar  
nonlinear Schr\"odinger equation \cite{LMSW96}. The eigenvalues of the 
saddle are: a positive one of order $\sqrt{\e}$, a negative one of order $\sqrt{\e}$, and a complex conjugate pair for which the real part is negative, and both real and imaginary parts are of order $\e$. The two eigenvalues of order $\sqrt{\e}$, correspond to a fish-shape dynamics, 
and the other two eigenvalues correspond to a stable spiraling dynamics.
\item The dynamics in the invariant diagonal subspace given by $p=q$, is 
governed by the PNLS. The unstable manifold of the saddle resides in the 
invariant subspace. This leads to that the persistent homoclinic orbit 
for the PNLS is the only one for the PVNLS, asymptotic to the saddle.
\end{enumerate}

\newpage
\eqnsection{Dynamics in the 4D Invariant Subspace of Functions Independent 
of $x$}

The 4D subspace 
\begin{equation}
\Pi = \{ (p,q)\ | \ \pa_x p = \pa_x q = 0 \} \ , 
\label{Pi}
\end{equation}
is invariant under the PVNLS flow (\ref{pnls1},\ref{pnls2}). Dynamics in $\Pi$ is governed by 
\begin{eqnarray}
& & ip_t + \frac{1}{2} [(|p|^2 + |q|^2)-\om^2] p = 
i \e [ -\al p - \be ]\ , \label{ppi1}\\
& & iq_t + \frac{1}{2} [(|p|^2 + |q|^2)-\om^2] q = 
i \e [ -\al q - \be ]\ , \label{ppi2}
\end{eqnarray}

\subsection{Unperturbed Dynamics}

When $\e =0$, the unperturbed dynamics is governed by
\begin{eqnarray}
& & ip_t + \frac{1}{2} [(|p|^2 + |q|^2)-\om^2] p = 0\ , \label{upi1}\\
& & iq_t + \frac{1}{2} [(|p|^2 + |q|^2)-\om^2] q = 0\ , \label{upi2}
\end{eqnarray}
Introducing the polar coordinates
\begin{equation}
p = \sqrt{I_1} e^{i \th_1 }\ , \ \ q = \sqrt{I_2} e^{i \th_2 }\ ,
\label{polar}
\end{equation}
one gets the expression for the general solution,
\begin{eqnarray*}
& & I_1 = \ \mbox{constant}\ , \ \ I_2 = \ \mbox{constant}\ , \\
& & \th_1 = \frac{1}{2}(I_1 + I_2 -\om^2)t + \ga_1\ , \\
& & \th_2 = \frac{1}{2}(I_1 + I_2 -\om^2)t + \ga_2\ .
\end{eqnarray*}
Equilibria consist of the origin $p=q=0$ and the submanifold              
\begin{equation}
\E = \{ (I_1,I_2,\th_1,\th_2) \ | \ I_1 + I_2 = \om^2 \} \ .
\label{eman}
\end{equation}
For an illustration of $\E$, see Figure \ref{efig}. 
\begin{figure}
\vspace{1.5in}
\caption{An illustration of the equilibrium manifold $\E$.}
\label{efig}
\end{figure}

\subsection{Perturbed Dynamics}

In terms of the polar coordinates (\ref{polar}), Equations (\ref{ppi1},\ref{ppi2}) can be rewritten as 
\begin{eqnarray}
\dot{I}_1 &=& \e [ -2 \al I_1 - 2 \be \sqrt{I_1} \cos \th_1 ] \ , 
\label{pol1} \\
\dot{\th}_1 &=& \frac{1}{2} [I_1 + I_2 - \om^2 ] +\e \be \frac{\sin \th_1} 
{\sqrt{I_1}}\ , 
\label{pol2} \\
\dot{I}_2 &=& \e [ -2 \al I_2 - 2 \be \sqrt{I_2} \cos \th_2 ] \ , 
\label{pol3} \\
\dot{\th}_2 &=& \frac{1}{2} [I_1 + I_2 - \om^2 ] +\e \be \frac{\sin \th_2} 
{\sqrt{I_2}}\ . \label{pol4}
\end{eqnarray}
\begin{lemma}
All the fixed points of these equations (\ref{pol1}-\ref{pol4}) 
satisfy the conditions
\[
I_1 = I_2 \ , \ \ \th_1 = \th_2 \ .
\]
\label{ple1}
\end{lemma}
Proof: Setting the right hand side of (\ref{pol1}-\ref{pol4}), one 
can deduce that for fixed points, 
\[
\tan \th_1 = \tan \th_2 \ , \ \ \cos \th_1 \cos \th_2 > 0 \ , 
\ \  \sin \th_1 \sin \th_2 > 0 \ .
\]
Thus $\th_1 = \th_2$, which implies $I_1 = I_2$. $\Box$

As in the case of perturbed scalar nonlinear Schr\"odinger equation 
\cite{LMSW96}, there are three fixed points, 
\begin{enumerate}
\item The focus $O_\e$ in the neighborhood of origin,
\begin{equation}
I_1 = \e^2 4 \be^2 \om^{-4} + \cdot \cdot \cdot \ , \ \ 
\cos \th_1 = - \frac{\al}{\be} \sqrt{I_1} \ , \ \ 
\be \sin \th_1 > 0\ .
\label{oe} 
\end{equation}
The eigenvalues are 
\begin{eqnarray}
\mu_{1,2} &=& - \e \al \pm i \sqrt{\left ( \frac {\e \be \sin \th_1 }
{\sqrt{I_1}}\right )^2 - 2 \e \be \sqrt{I_1} \sin \th_1 } \ , 
\label{evoe1} \\
\mu_{3,4} &=& - \e \al \pm i \e \frac {\be \sin \th_1 }
{\sqrt{I_1}} \ , \label{evoe2}
\end{eqnarray}
where $I_1$ and $\th_1$ are given in (\ref{oe}).
\item The focus $P_\e$ in the neighborhood of $I_1 = I_2 = \frac{1}{2} \om^2$,
\begin{equation}
I_1 = \frac{1}{2} \om^2 +\e \om^{-1} \sqrt{2 \be^2 - \al^2 \om^2}
+ \cdot \cdot \cdot \ , \ \ 
\cos \th_1 = - \frac{\al}{\be} \sqrt{I_1} \ , \ \ 
\be \sin \th_1 < 0\ .
\label{pe} 
\end{equation}
The eigenvalues are 
\begin{eqnarray}
\mu_{1,2} &=& - \e \al \pm i \sqrt{\e}\sqrt{- 2 \be \sqrt{I_1} \sin \th_1
+\e \left ( \frac {\be \sin \th_1 }{\sqrt{I_1}} \right )^2 }
\label{evpe1} \\
\mu_{3,4} &=& - \e \al \pm i \e \frac {\be \sin \th_1 }
{\sqrt{I_1}} \ , \label{evpe2}
\end{eqnarray}
where $I_1$ and $\th_1$ are given in (\ref{pe}).
\item The saddle $Q_\e$ in the neighborhood of $I_1 = I_2 = \frac{1}{2} \om^2$,
\begin{equation}
I_1 = \frac{1}{2} \om^2 -\e \om^{-1} \sqrt{2 \be^2 - \al^2 \om^2}
+ \cdot \cdot \cdot \ , \ \ 
\cos \th_1 = - \frac{\al}{\be} \sqrt{I_1} \ , \ \ 
\be \sin \th_1 > 0\ .
\label{qe} 
\end{equation}
The eigenvalues are 
\begin{eqnarray}
\mu_{1,2} &=& - \e \al \pm \sqrt{\e}\sqrt{2 \be \sqrt{I_1} \sin \th_1
-\e \left ( \frac {\be \sin \th_1 }{\sqrt{I_1}} \right )^2 }
\label{evqe1} \\
\mu_{3,4} &=& - \e \al \pm i \e \frac {\be \sin \th_1 }
{\sqrt{I_1}} \ , \label{evqe2}
\end{eqnarray}
where $I_1$ and $\th_1$ are given in (\ref{qe}).
\end{enumerate}

The pair of homoclinic orbits to be located, are asymptotic to the 
saddle $Q_\e$. Inside the subspace $\Pi$, 
attention will be focused upon the neighborhood of $I_1 = I_2 = 
\frac{1}{2} \om^2$.
Then naturally one wants to introduce coordinates ($J_1,J_2$) in this neighborhood,
\begin{equation}
I_1 =  \frac{1}{2}\om^2 + J_1\ , \ \ I_2 = \frac{1}{2} \om^2 + J_2\ ,
\label{coJ}
\end{equation}
and Equations (\ref{pol1}-\ref{pol4}) can be rewritten as
\begin{eqnarray}
\dot{J}_1 &=& \e \left [ - \al \om^2 -2 \al J_1 - 2 \be 
\sqrt{\frac{1}{2} \om^2 + J_1} \cos \th_1\right ] \ , 
\label{Jco1} \\
\dot{\th}_1 &=&\frac{1}{2}  (J_1 + J_2) +\e \be \frac{\sin \th_1} 
{\sqrt{\frac{1}{2} \om^2 + J_1}}\ , 
\label{Jco2} \\
\dot{J}_2 &=& \e \left [ - \al \om^2 -2 \al J_2 - 2 \be 
\sqrt{ \frac{1}{2}\om^2 + J_2} \cos \th_2 \right ] \ , 
\label{Jco3} \\
\dot{\th}_2 &=& \frac{1}{2} (J_1 + J_2) +\e \be \frac{\sin \th_2} 
{\sqrt{\frac{1}{2} \om^2 + J_2}}\ . 
\label{Jco4} 
\end{eqnarray}
Rescaling the variables as follows, one can reveal the dynamics at different scales better. Let
\begin{equation}
J_1 = \sqrt{\e} j_1\ , \ \ J_2 = \sqrt{\e} j_2\ ,
\ \ \tau = \sqrt{\e} t \ ,
\label{colj}
\end{equation}
and Equations (\ref{Jco1}-\ref{Jco4}) can be rewritten as
\begin{eqnarray}
j'_1 &=& - \al \om^2 -\sqrt{2} \be \om \cos \th_1 - \sqrt{\e} 
2 \al j_1 - 2 \be \left (\sqrt{\frac{1}{2} \om^2 + \sqrt{\e}j_1} - 
\frac {\om}{\sqrt{2}} \right )\cos \th_1  \ , 
\label{ljco1} \\
\th'_1 &=& \frac{1}{2} (j_1 + j_2) + \sqrt{\e} \be \frac{\sin \th_1} 
{\sqrt{\frac{1}{2} \om^2 + \sqrt{\e}j_1}}\ , 
\label{ljco2} \\
j'_2 &=& - \al \om^2 -\sqrt{2} \be \om \cos \th_2 - \sqrt{\e} 
2 \al j_2 - 2 \be \left (\sqrt{\frac{1}{2} \om^2 + \sqrt{\e}j_2} - \frac {\om}{\sqrt{2}} \right )\cos \th_2  \ , 
\label{ljco3} \\
\th'_2 &=&  \frac{1}{2}(j_1 + j_2) + \sqrt{\e} \be \frac{\sin \th_2} 
{\sqrt{ \frac{1}{2}\om^2 + \sqrt{\e}j_2}}\ , 
\label{ljco4} 
\end{eqnarray}
where $' = \frac {d} {d \tau}$. 

\subsubsection{Leading Order Dynamics} 

Only keeping the leading order terms on the right hand side of (\ref{ljco1}-\ref{ljco4}), one gets 
\begin{eqnarray}
j'_1 &=& - \al \om^2 -\sqrt{2} \be \om \cos \th_1 \ ,
\label{lead1} \\
\th'_1 &=& \frac{1}{2} (j_1 + j_2) \ ,
\label{lead2} \\
j'_2 &=& - \al \om^2 -\sqrt{2} \be \om \cos \th_2 \ ,
\label{lead3} \\
\th'_2 &=&  \frac{1}{2}(j_1 + j_2) \ .
\label{lead4} 
\end{eqnarray}
This system is an integrable Hamiltonian system with the 
Hamiltonian
\begin{equation}
\H = \frac{1}{4}(j_1 +j_2)^2 +\al \om^2 (\th_1 +\th_2) +
\sqrt{2} \be \om (\sin \th_1 +\sin \th_2 )\ ,
\label{leadH}
\end{equation} 
and another functionally independent constant of motion,
\begin{equation}
K = \th_1 - \th_2 \ .
\label{leadK}
\end{equation}
The Hamiltonian form of (\ref{lead1}-\ref{lead4}) is
\[
j'_1 = - \frac {\pa \H}{\pa \th_1}\ , \ \ 
\th'_1 = \frac {\pa \H}{\pa j_1}\ , \ \ 
j'_2 = - \frac {\pa \H}{\pa \th_2}\ , \ \ 
\th'_2 = \frac {\pa \H}{\pa j_2}\ .
\]
To this order, four lines in the equilibrium manifold $\E$ (\ref{eman}) 
persist. These four lines of fixed points of (\ref{lead1}-\ref{lead4})
are
\begin{equation}
j_1 = - j_2 \ , \ \ \th_1 = \pm \arccos \bigg [ - \frac {\al \om}{\sqrt{2} \be } \bigg ] \ , \ \ \th_2 = \pm \arccos \bigg [ 
- \frac {\al \om}{\sqrt{2} \be } \bigg ] \ .
\label{4line}
\end{equation}
Two lines of fixed points for which $\th_1 \neq \th_2$ will disappear 
under the full perturbed flow (\ref{pol1}-\ref{pol4}). One point on each 
of the other two lines will persist and becomes $P_\e$ (\ref{pe}) or 
$Q_\e$ (\ref{qe}). The eigenvalues of the four lines of fixed points are
\begin{equation}
\mu_{1,2} = \pm \sqrt{ \frac {1} {\sqrt{2}} \be \om (\sin \th_1 + 
\sin \th_2 )} \ , \ \ \mu_{3,4} = 0 \ .
\label{evl}
\end{equation}
Denote by $l_u$ the line in (\ref{4line}), for which $\th_1 = \th_2$ and 
$\be \sin \th_1 >0$. For $l_u$, $\mu_{1,2}$ in (\ref{evl}) are real. 
Denote by $l_c$ the line in (\ref{4line}), for which $\th_1 = \th_2$ and 
$\be \sin \th_1 <0$. For $l_c$, $\mu_{1,2}$ in (\ref{evl}) are purely 
imaginary. Denote by $l_0^1$ and $l_0^2$ the rest two lines in (\ref{4line}), for which $\th_1 = -\th_2$. For $l_0^1$ and $l_0^2$, 
$\mu_{1,2}$ in (\ref{evl}) are zero. Denote by $Q_0$ the point on 
$l_u$ such that $j_1 = - j_2 = 0$. $Q_0$ has the coordinates
\begin{equation}
j_1 = - j_2 = 0 \ , \ \ 
\th_1 = \th_2 = \left \{ \begin{array}{l} 
\pi - \arctan \bigg [ \frac {\sqrt{2 \be^2 - \al^2 \om^2}} 
{\al \om }\bigg ] \ , \ \ \mbox{if} \ \be > 0 \ , \cr 
-\arctan \bigg [ \frac {\sqrt{2 \be^2 - \al^2 \om^2}} 
{\al \om }\bigg ] \ , \ \ \mbox{if} \ \be < 0 \ . \cr 
\end{array} \right .
\label{Q0}
\end{equation}
Denote by $P_0$ the point on 
$l_c$ such that $j_1 = - j_2 = 0$. $P_0$ has the coordinates
\begin{equation}
j_1 = - j_2 = 0 \ , \ \ 
\th_1 = \th_2 = \left \{ \begin{array}{l} 
-\pi + \arctan \bigg [ \frac {\sqrt{2 \be^2 - \al^2 \om^2}} 
{\al \om }\bigg ] \ , \ \ \mbox{if} \ \be > 0 \ , \cr 
\arctan \bigg [ \frac {\sqrt{2 \be^2 - \al^2 \om^2}} 
{\al \om }\bigg ] \ , \ \ \mbox{if} \ \be < 0 \ . \cr 
\end{array} \right .
\label{P0}
\end{equation}
\nid
Next we study the dynamics of (\ref{lead1}-\ref{lead4}). Denote 
by $\Pi_\Dl$ the 3D subspace of $\Pi$ (\ref{Pi}),
\begin{equation}
\Pi_\Dl = \{ (j_1,j_2,\th_1,\th_2) \ | \ 
K = \th_1 - \th_2 = \Dl \} \ .
\label{Chi}
\end{equation}
$\Pi_\Dl$  is an invariant subspace under the flow 
(\ref{lead1}-\ref{lead4}). Restricted to $\Pi_\Dl$, $\H$ (\ref{leadH}) 
takes the form similar to that in \cite{LMSW96},
\[
\H = \frac{1}{4} j^2 + 2 \al \om^2 \th + 2 \sqrt{2} \be \om \cos (\Dl /2) 
\sin \th \ , 
\]
where $j = j_1 + j_2$ and $\th = \th_1 - \Dl /2$.
Denote by $W^s_\Pi(l_u)$, $W^u_\Pi(l_u)$, and $W^c_\Pi(l_u)$ the 
2D stable, unstable, and center manifolds of $l_u$ respectively. 
$W^s_\Pi(l_u)= W^u_\Pi(l_u)$ has the topology shown in Figure \ref{fish}. $W^c_\Pi(l_u)$ has the topology shown in Figure \ref{clu}. Similar to $W^c_\Pi(l_u)$, there is also an invariant 
2D submanifold $W^c_\Pi(l_c)$ around $l_c$ as shown in Figure \ref{clc}. $W^s_\Pi(l_u)= W^u_\Pi(l_u)$ is given by the singular 
level set of $\H$ in $\Pi_{\Dl = 0}$. In $\Pi_\Dl$  ($\Dl \neq 0$), similar 
picture holds as shown in Figure \ref{nfish}, but with drifting 
along the asymptotic line. Denote by $W^{cs}_\Pi(l_u)$ and $W^{cu}_\Pi(l_u)$ the 3D center-stable and center-unstable manifolds of $l_u$. $W^{cs}_\Pi(l_u) = W^{cu}_\Pi(l_u)$ has the topology shown in Figure \ref{csu}.
\begin{figure}
\vspace{1.5in}
\caption{The stable and unstable manifolds of $l_u$.}
\label{fish}
\end{figure}
\begin{figure}
\vspace{1.5in}
\caption{The center manifold of $l_u$.}
\label{clu}
\end{figure}
\begin{figure}
\vspace{1.5in}
\caption{The center manifold of $l_c$.}
\label{clc}
\end{figure}
\begin{figure}
\vspace{1.5in}
\caption{Stable and unstable manifolds in $\Pi_\Dl$  ($\Dl \neq 0$).}
\label{nfish}
\end{figure}
\begin{figure}
\vspace{1.5in}
\caption{The center-stable and center-unstable manifolds of $l_u$.}
\label{csu}
\end{figure}
\begin{figure}
\vspace{1.5in}
\caption{Persistent center-stable and center-unstable manifolds.}
\label{pcsu}
\end{figure}

\subsubsection{Full Dynamics}

Since $W^c_\Pi(l_u)$ is normally hyperbolic under the leading order 
flow (\ref{lead1}-\ref{lead4}), it persists under the full perturbed 
flow (\ref{ljco1}-\ref{ljco4}). Denote this persistent 2D locally invariant center manifold by $W^{(c,\e)}_\Pi(l_u)$. Denote the 
persistent 3D locally invariant center-stable and center-unstable manifolds by $W^{(cs,\e)}_\Pi(l_u)$ and $W^{(cu,\e)}_\Pi(l_u)$.
$W^{(cs,\e)}_\Pi(l_u)$ and $W^{(cu,\e)}_\Pi(l_u)$ have the topology 
shown in Figure \ref{pcsu}. Notice that $W^{(cs,\e)}_\Pi(l_u)$ 
forms the stable manifold of the saddle $Q_\e$ (\ref{qe}), which we also denote by $W^s_\Pi(Q_\e)$. The unstable manifold of $Q_\e$, denoted by 
$W^u_\Pi(Q_\e)$, is one dimensional.

\newpage
\eqnsection{Global Integrable Theory}

Consider the integrable vector nonlinear Schr\"odinger equations
($\e =0$ in (\ref{pnls1},\ref{pnls2})),
\begin{eqnarray}
& & ip_t + p_{xx} + \frac{1}{2} [(|p|^2 + |q|^2)-\om^2] p = 0\ , \label{nls1}\\
& & iq_t + q_{xx} + \frac{1}{2} [(|p|^2 + |q|^2)-\om^2] q = 0\ . \label{nls2}
\end{eqnarray}
The Lax pair is given as \cite{Man74},
\begin{equation}
\psi_x = U \psi \ , \ \ \psi_t = V \psi \ ,
\label{lax}
\end{equation}
where
\[
U = \la A_0 + A_1 \ , \ \ V = \la^2 A_0 + \la A_1 + A_2 \ , 
\]
\[
A_0 = \left ( \begin{array}{lcr} -\frac {2}{3} i & 0 & 0 \cr 
0& \frac {1}{3} i & 0 \cr 
0& 0 &\frac {1}{3} i \cr \end{array} \right ) \ , \ \ 
A_1 = \left ( \begin{array}{lcr} 0 & \frac {1}{2} p & \frac {1}{2} q \cr 
-\frac {1}{2} \bp & 0 & 0 \cr 
-\frac {1}{2} \bq & 0 & 0 \cr 
\end{array} \right ) \ , 
\]
\[
A_2 = - \frac {i}{4}\left ( \begin{array}{lcr}
-(|p|^2+|q|^2)+4\om^2 & - 2p_x & -2 q_x \cr 
-2 \bp_x & |p|^2+2\om^2 & q\bp \cr 
-2\bq_x & \bq p & |q|^2+2 \om^2\cr 
\end{array} \right ) \ .
\]

\subsection{Linearized Equations}

We start with the solution of (\ref{upi1},\ref{upi2}),
\begin{equation}
p = a e^{i \dl_1} \ , \ \ q = b e^{i \dl_2} \ ,
\label{ups}
\end{equation}
where
\[
\dl_1 = \frac{1}{2} (a^2 +b^2 -\om^2) t + \ga_1 \ , \ \ 
\dl_2 = \frac{1}{2} (a^2 +b^2 -\om^2) t + \ga_2 \ .
\]
Linearizing (\ref{nls1},\ref{nls2}) at this solution, by setting 
\[
p = e^{i \dl_1}[a+\tp ]\ , \ \ q = e^{i \dl_2}[b + \tq ]\ ,
\]
one gets
\begin{eqnarray*}
& & i\tp_t + \tp_{xx} + \frac{1}{2} a[ a(\tp + \bar{\tp})+ 
b (\tq + \bar{\tq})]  = 0\ , \\
& & i\tq_t + \tq_{xx} + \frac{1}{2} b[ a(\tp + \bar{\tp})+ 
b (\tq + \bar{\tq})] = 0\ . 
\end{eqnarray*}
Setting 
\[
\tp = f_+ e^{ikx+ \Om t} + \overline{f_-} e^{-ikx+ \bar{\Om}t} \ , 
\ \ 
\tq = g_+ e^{ikx+ \Om t} + \overline{g_-} e^{-ikx+ \bar{\Om}t} \ ,
\ \ (k \in Z) 
\]
one gets the dispersion relations
\begin{equation}
\Om_{1,2} = \pm k \sqrt{a^2 + b^2 - k^2}\ , \ \ 
\Om_{3,4} = \pm ik^2 \ .
\label{grate}
\end{equation}
Thus, if we choose $1 < \sqrt{a^2+b^2} < 2$, then there is only one 
unstable mode $\cos x$. This is the reason that we restricts $\om$ as 
in (\ref{inls1},\ref{inls2}). When $k=1$, 
\begin{equation}
\Om_{1,2} = \pm \sg \ , \ \ \sg = \sqrt{ a^2 + b^2 - 1}\ .
\label{luc}
\end{equation}
The corresponding eigenfunctions are
\[
\tp = a e^{\pm i \vphi }e^{\pm \sg t} \cos x \ ,
\ \ \tq = b e^{\pm i \vphi }e^{\pm \sg t} \cos x \ ,
\ \ \vphi = \arctan \sg \ .
\]

\subsection{Homoclinic Orbits and Figure Eight Structures}

The B\"acklund-Darboux transformation given in 
\cite{WF00} \cite{FSW00} will be utilized to construct homoclinic 
orbits asymptotic to the periodic orbits (\ref{ups}) up to phase shifts.
\begin{theorem}[\cite{WF00} \cite{FSW00}]
Let $p$ and $q$ be a solution of the vector nonlinear Schr\"odinger 
equations (\ref{nls1},\ref{nls2}), and let $\psi = 
(\psi_1,\psi_2,\psi_3)^T$ 
be an eigenfunction 
of the Lax pair (\ref{lax}) at ($\la, p,q$). If one defines 
\begin{eqnarray}
\hp &=& p + 2 i (\bar{\la} - \la ) \frac {\psi_1 \overline{\psi_2}}
{|\psi_1|^2 + |\psi_2|^2 + |\psi_3|^2 } \ ,  \label{bd1}\\ 
\hq &=& q + 2 i (\bar{\la} - \la ) \frac {\psi_1 \overline{\psi_3}}
{|\psi_1|^2 + |\psi_2|^2 + |\psi_3|^2 } \ , \label{bd2}
\end{eqnarray}
then $\hp$ and $\hq$ also solve the vector nonlinear Schr\"odinger 
equations (\ref{nls1},\ref{nls2}).
\label{bdthm}
\end{theorem}
We solve the Lax pair (\ref{lax}) at ($\la, p,q$) with $p$ and $q$ given by the periodic orbits (\ref{ups}), we get
\begin{equation}
\psi^{\pm} = \left ( \begin{array}{c} \psi_1^{\pm}\cr
 \psi_2^{\pm} \cr \psi_3^{\pm} \cr \end{array} \right ) 
= \left ( \begin{array}{c} (i\la \mp i)e^{i(\dl_1+\dl_2)} \cr 
a e^{i\dl_2} \cr b e^{i\dl_1} \cr \end{array}\right ) 
e^{\pm \sg_1 t +i \sg_2 t} e^{\pm i \k^{\pm} x} \ ,
\label{efunc}
\end{equation}
where
\begin{eqnarray}
& & \k^{\pm} = - \frac{1}{6} \la \pm \frac{1}{2} 
\sqrt{a^2+b^2 + \la^2}\ , \ \ \sqrt{a^2+b^2 + \la^2} = 1\ , 
\label{wno} \\
& & \sg_1 = \frac{i\la}{2} \ , \ \ 
\sg_2 = - \frac{1}{12} [7(a^2+b^2)+2] \ . 
\label{freq} 
\end{eqnarray}
\begin{remark}
From (\ref{freq}), we see that in order to have temporal growth, 
the imaginary part of $\la$ can not be zero. Also in order to have 
a nontrivial B\"acklund-Darboux transformation (Theorem \ref{bdthm}),
the imaginary part of $\la$ can not be zero. From (\ref{wno}), if 
$a^2+b^2 >1$, then $\la$ is purely imaginary. If we also require 
$a^2+b^2 <4$, i.e., in the one unstable mode regime (\ref{grate}), 
then $\sqrt{a^2+b^2 + \la^2}\geq 2$ will lead to real $\la$ and no 
temporal growth.
\end{remark}
\begin{remark}
Notice that the imaginary part of $\k^{\pm}$ is not zero, thus 
$\psi^{\pm}$ are not periodic or antiperiodic functions in $x$. 
In fact, they grow or decay in $x$. This fact shows that complex 
double points are not necessary in constructing homoclinic orbits 
through B\"acklund-Darboux transformations (see also \cite{Li00a}).
\end{remark}
Let 
\begin{equation}
\psi = c^+ \psi^+ + c^- \psi^- \ , 
\label{gfunc}
\end{equation}
where $c^+$ and $c^-$ are two arbitrary complex constants.
We introduce the notations
\begin{eqnarray}
& & c^+/c^- = e^{\rho + i \vth }\ , \ \ \la = -i \sg \ , \ \
\sg = \sqrt{a^2+b^2 - 1}\ , \label{nota1} \\
& & \sg + i = \sqrt{a^2+b^2} e^{i (\frac{\pi}{2} - \vth_0)} \ , 
\ \ \vth_0 = \arctan \sg \ , \ \
\sg_1 = \sg /2 \ , \label{nota2}\\
& & \tau = \frac{\sg}{2} t + \frac{\rho}{2}\ , \ \
y = \frac{1}{2} x +\frac{\vth}{2}\ , \label{nota3}
\end{eqnarray}
where $\rho$ and $\vth$ are called the B\"acklund parameters. 
If one chooses $\la = i \sg$, one will get the same result.
We can rewrite $\psi$ in the following form,
\begin{eqnarray*}
\psi_1 &=& \bigg [ \cosh \tau \cos (y - \frac{\pi}{2} + \vth_0) + i \sinh \tau 
\sin (y - \frac{\pi}{2} + \vth_0) \bigg ] 2 \sqrt{c^+c^-} \sqrt{a^2+b^2} 
e^{i(\dl_1+\dl_2)}e^{i\sg_2 t} e^{-\frac{1}{6}\sg x}\ , \\
\psi_2 &=& \bigg [ \cosh \tau \cos y  + i \sinh \tau 
\sin y \bigg ] 2 \sqrt{c^+c^-} a 
e^{i\dl_2}e^{i\sg_2 t} e^{-\frac{1}{6}\sg x}\ , \\
\psi_3 &=& \bigg [ \cosh \tau \cos y  + i \sinh \tau 
\sin y \bigg ] 2 \sqrt{c^+c^-} b 
e^{i\dl_1}e^{i\sg_2 t} e^{-\frac{1}{6}\sg x}\ ,
\end{eqnarray*}
where the factor $2 \sqrt{c^+c^-} e^{i\sg_2 t} e^{-\frac{1}{6}\sg x}$ 
will be canceled away. By (\ref{bd1},\ref{bd2}), one gets
\begin{equation}
\hp = a e^{i\dl_1} h\ , \ \ \hq = b e^{i\dl_2} h\ , 
\label{horb}
\end{equation}
where
\begin{eqnarray}
h &=& \bigg [ \cos (2\vth_0) + i \sin (2\vth_0) \tanh (2\tau) -
\sin \vth_0 \ \mbox{sech} (2\tau) \cos (2y + \vth_0 - \frac{\pi}{2})
\bigg ] \nonumber \\
& & \bigg [ 1 + \sin \vth_0 \ \mbox{sech} (2\tau) \cos (2y + 
\vth_0 - \frac{\pi}{2}) \bigg ]^{-1} \ ,
\label{exh}
\end{eqnarray}
where $\dl_1$ and $\dl_2$ are given in (\ref{ups}), other notations are 
given in (\ref{nota1}-\ref{nota3}). For the moment, the 
even-in-$x$ condition has not been inforced. For fixed $a$ and $b$,
the phases $\ga_1$ and $\ga_2$ (\ref{ups}), and the B\"acklund 
parameters $\rho$ and $\vth$ parametrize a four dimensional submanifold 
with a figure eight structure. If $a=b$ and $\ga_1=\ga_2$, then the homoclinic 
orbit reduces to that for the scalar nonlinear Schr\"odinger equation 
\cite{Li01a}.

As $t \ra \pm \infty$, 
\begin{equation}
h \ra e^{\pm i 2 \vth_0 } \ .
\label{symph}
\end{equation}
The even-in-$x$ condition can be achieved by restricting the B\"acklund 
parameter $\vth$ to special values. That is, if one requires that 
\begin{equation}
\vth + \vth_0 - \frac{\pi}{2} = 0 \ , \ \pi \ , 
\label{evc}
\end{equation}
then $h$ is even in $x$, 
\begin{eqnarray}
h &=& \bigg [ \cos (2\vth_0) + i \sin (2\vth_0) \tanh (2\tau) \mp
\sin \vth_0 \ \mbox{sech} (2\tau) \cos x
\bigg ] \nonumber \\
& & \bigg [ 1 \pm \sin \vth_0 \ \mbox{sech} (2\tau) \cos x
 \bigg ]^{-1} \ ,
\label{exeh}
\end{eqnarray}
where the upper sign corresponds to $\vth + \vth_0 - \frac{\pi}{2}= 0$.
For fixed $a$ and $b$,
the phases $\ga_1$ and $\ga_2$ (\ref{ups}), and the B\"acklund 
parameter $\rho$ parametrize a three dimensional submanifold 
with a figure eight structure. If one further chooses $\ga_1 = \ga_2$ (\ref{ups}), then (\ref{ups}) represents a periodic orbit. In such case,
the phase $\ga_1 = \ga_2$ and the B\"acklund 
parameter $\rho$ parametrize a two dimensional submanifold 
with a figure eight structure. 

\subsection{A Melnikov Vector}

The Hamiltonian for the integrable vector nonlinear Schr\"odinger equations (\ref{nls1},\ref{nls2}) is given as,
\begin{equation}
H = \int^{2\pi}_0 \{ |p_x|^2+|q_x|^2 - \frac{1}{4} [ (|p|^2+|q|^2)^2 
-2\om^2 (|p|^2+|q|^2)]\} dx\ .
\label{hami}
\end{equation}
(\ref{nls1},\ref{nls2}) can be rewritten in the Hamiltonian form, 
\[
ip_t = \frac{\dl H}{\dl \bp}\ , \ \ 
i\bp_t = -\frac{\dl H}{\dl p}\ , \ \ 
iq_t = \frac{\dl H}{\dl \bq}\ , \ \ 
i\bq_t = -\frac{\dl H}{\dl q}\ .
\]
The two $L^2$ norms 
\begin{equation}
E_1 = \int^{2\pi}_0 |p|^2 dx\ , \ \ 
E_2 = \int^{2\pi}_0 |q|^2 dx\ ,
\label{2nor}
\end{equation}
are also invariant functionals. We will build the Melnikov vector 
by the gradient of the invariant functional $G$ obtained through 
a combination out of $H$, $E_1$, and $E_2$,
\begin{eqnarray}
G &=& H +\frac{1}{2} \left [ \frac{1}{4\pi }(E_1+E_2) - \om^2 \right ]
(E_1+E_2) \nonumber \\
&=& \int^{2\pi}_0 \{ |p_x|^2+|q_x|^2 - \frac{1}{4} (|p|^2+|q|^2)^2 
\} dx + \frac{1}{8\pi }(E_1+E_2)^2\ . \label{mff}
\end{eqnarray}

\newpage
\eqnsection{Persistent Homoclinic Orbits}

\subsection{The Invariant Diagonal Subspace}

The diagonal subspace
\[
\Sg = \{ (p,q)\ | \ p=q \}
\]
is invariant under the PVNLS flow (\ref{pnls1})-(\ref{pnls2}). The 
dynamics in $\Sg$ is governed by the PNLS
\begin{equation}
iq_t + q_{xx} + \frac{1}{2} [2|q|^2-\om^2] q = 
i \e [ q_{xx} -\al q - \be ]\ , 
\label{pnls}
\end{equation}
studied in \cite{Li01a}.

\subsection{The Main Theorem}

\begin{theorem}[Main Theorem]
There exists a $\e_0 > 0$, such that 
for any $\e \in (0, \e_0)$, there exists 
a codimension 1 surface in the space of $(\alpha,\beta, \om) \in 
R^+\times R^+\times R^+$ where 
$\om \in (1, 2)/S$, $S$ is a finite subset, and 
$\al \om < \sqrt{2} \be$. For any $(\alpha ,\beta, \omega)$ on the 
codimension 1 surface, the singularly perturbed vector nonlinear 
Schr\"odinger equations (\ref{pnls1})-(\ref{pnls2}) possesses a 
homoclinic orbit asymptotic to the saddle
$Q_\epsilon$ (\ref{qe}). This orbit is also the homoclinic orbit 
for the singularly perturbed scalar nonlinear Schr\"odinger equation
(\ref{pnls}) studied in \cite{Li01a}, and is the only one asymptotic 
to the saddle $Q_\epsilon$ (\ref{qe}) for the singularly perturbed 
vector nonlinear Schr\"odinger equations (\ref{pnls1})-(\ref{pnls2}).
The codimension 1 surface has the 
approximate representation $\al =\al (\om)$ given in Figure \ref{alcur}, 
obtained in \cite{Li01a}.
\end{theorem} 

Proof. Existence of a homoclinic orbit has been proved in \cite{Li01a}.
Under the PVNLS flow (\ref{pnls1})-(\ref{pnls2}), the saddle $Q_\epsilon$ 
(\ref{qe}) has a two dimensional unstable manifold $W^u(Q_\e)$ which is 
included in the invariant diagonal subspace $\Sg$,
\[
W^u(Q_\e) \subset \Sg \ .
\]
Thus any homoclinic orbit asymptotic to $Q_\epsilon$ has to be inside 
$\Sg$. $\Box$

\begin{figure}
\vspace{1.5in}
\caption{The graph of $\al = \al (\om)$.}
\label{alcur}
\end{figure}

\begin{remark}
Following the same arguments as in \cite{Li01a}, Unstable Fiber Theorem 
and Center-Stable Manifold Theorem can be easily established for the 
PVNLS (\ref{pnls1})-(\ref{pnls2}).
\end{remark}

\subsection{An Alternative Melnikov Function}

In \cite{Li01a}, the Melnikov function was built with the Floquet 
discriminant. Here we give an alternative Melnikov function built 
with the invariant functional $G$:
\[
M= \int_{-\infty}^{+\infty} \frac {dG} {dt} dt = \int_{-\infty}^{+\infty} 
\int_0^{2\pi} \left [ \frac{\dl G}{\dl p}p_t + \frac{\dl G}{\dl \bar{p}}\bar{p}_t 
+\frac{\dl G}{\dl q}q_t +\frac{\dl G}{\dl \bar{q}}\bar{q}_t \right ]dx dt \ ,
\]
where the integral is evaluated along the homoclinic orbit 
(\ref{horb},\ref{exeh})
at $a = b =\frac{1}{\sqrt{2}} \om$ and $\dl_1 = \ga_1 = \ga_2 = \dl_2$,
$\pm$ in (\ref{exeh}) leads to the same result.
\begin{equation}
M = \e \left (\frac{\sg}{2}\right )^{-1} \bigg [ M_1 - \al M_2 - 
\be 2 \sqrt{2} 
\om \cos \ga_1 M_3 \bigg ] \ ,
\label{mel1}
\end{equation}
where 
\[
M_1 = \frac{1}{2} \om^2 \int_{-\infty}^{+\infty} \int_0^{2\pi} 
\bigg [ \om^2 (1-|h|^2)(\bh h_{xx} + h \bh_{xx})- 4 |h_{xx}|^2 
\bigg ] dx d\tau \ ,
\]
\[
M_2 = \om^2 \int_{-\infty}^{+\infty} \int_0^{2\pi} 
\bigg [ \om^2 |h|^2(1-|h|^2) - (\bh h_{xx}+h \bh_{xx}) \bigg ] 
dx d\tau \ ,
\]
\[
M_3 = \int_{-\infty}^{+\infty} \int_0^{2\pi} 
\bigg [ \frac{1}{2} \om^2 h^{(r)} (1-|h|^2) - h^{(r)}_{xx} 
\bigg ] dx d\tau \ ,
\]
where $h^{(r)}$ is the real part of $h$ (\ref{exeh}).

In \cite{Li01a}, we have showed that $M=0$ together with 
\begin{equation}
\H(0,0,\ga_1 + 2 \vth_0) - \H(0,0,\ga_1 - 2 \vth_0) = 0 \ ,
\label{sec1}
\end{equation}
give the function $\al =\al (\om)$ which predicts the leading 
order location of the homoclinic orbit in the external parameter 
space. Here $\pm 2 \vth_0$ are the asymptotic phases (\ref{symph}),
and the Hamiltonian $\H$ is given in (\ref{leadH}), 
\[
\H(j_1,j_2,\th_1) = \frac{1}{4}(j_1+j_2)^2 +2\al \om^2 \th_1 
+2\sqrt{2}\be \om \sin \th_1 \ .
\]
(\ref{sec1}) leads to
\begin{equation}
\cos \ga_1 = -\frac{\sqrt{2} \al \om \vth_0}{\be \sin 2\vth_0} \ .
\label{sec2}
\end{equation}
Combining (\ref{sec2}) with $M = 0$, one gets 
\begin{equation}
\al = \al (\om) = M_1 \bigg [ M_2 -\frac{4\om^2 \vth_0}{\sin 2 \vth_0}
M_3 \bigg ]^{-1} \ ,
\label{param}
\end{equation}
which has the same graph as in Figure \ref{alcur}. Thus the Melnikov 
function above is equivalent to that in \cite{Li01a}.

\eqnsection{Conclusion and Discussion}

The instability of the singularly perturbed vector nonlinear 
Schr\"odinger equations (\ref{pnls1})-(\ref{pnls2}) happens along the 
diagonal direction ($p=q$), which leads to the fact that the homoclinic 
orbit in the invariant diagonal subspace ($p=q$) is the only one for 
(\ref{pnls1})-(\ref{pnls2}). Existence of a homoclinic orbit in the 
invariant diagonal subspace ($p=q$) governed by the singularly perturbed
scalar nonlinear Schr\"odinger equation (\ref{pnls}) has been established 
in \cite{Li01a}.

If one releases the even constraint $p(t, -x) = p(t, x)$ and 
$q(t, -x) = q(t, x)$, the homoclinic orbit is given by (\ref{exh}) in 
which the spatial variable $x$ appears in the combination:
\[
x + \vth + \vth_0 - \frac{\pi}{2} \ ,
\]
where $\vth$ is the B\"acklund phase-parameter. (Even constraint is realized 
by setting $\vth + \vth_0 - \frac{\pi}{2} = 0\ \ \mbox{or} \ \ \pi$.) 
$\vth$ characterizes the $x$-location of the top of the hump $h(t,x)$ for any 
$t$. It is clear that the Melnikov integrals evaluated along $h(t,x)$ are 
independent of $\vth$. In fact, since the PVNLS (\ref{pnls1},\ref{pnls2}) 
is translationally invariant in $x$, the persistent homoclinic orbit is not 
isolated, if it exists. In the current study, the difficulty associated with 
noneveness, that we are facing, is that we need two Melnikov integrals due to 
the extra phase-parameter $\vth$. Using the real and imaginary parts of the 
constant of motion $F_1$ (see (3.22) in \cite{Li01a}), we can build two 
Melnikov integrals. Unfortunately, we found numerically that the second 
Melnikov integral built with the imaginary part of $F_1$ is identically zero. 
We do not know whether or not one can build two non-identically-zero 
Melnikov integrals. If not, second or higher order Melnikov integral 
might be needed.

The matter is much more important than just a technical difficulty. 
In this noneven case, Ablowitz et al. \cite{AHS96} \cite{AHHS00} 
found a novel scenario in which chaos occurs, that is, chaotic flipping 
between left and right travelling waves. This scenario is generic 
with respect to perturbations \cite{AHS96} \cite{AHHS00}. When even 
constraint is not imposed, this new scenario of chaotic dynamics is 
observed in experiments \cite{AHHS00}. With the difficulty mentioned in 
the last paragraph and the new results of Ablowitz et al. \cite{AHS96} 
\cite{AHHS00}, understanding the chaotic dynamics in the noneven case 
becomes a much more challenging and significant problem.

\vspace{0.3in}
{\bf Acknowledgement:} The author is greatly indebted to Professor 
Mark Ablowitz and Professor Jianke Yang for fruitful discussions 
which led to the current study.

\vspace{0.2in}
{\bf Figure Captions:}

Figure 2.1: An illustration of the equilibrium manifold $\E$.

Figure 2.2: The stable and unstable manifolds of $l_u$.

Figure 2.3: The center manifold of $l_u$.

Figure 2.4: The center manifold of $l_c$.

Figure 2.5: Stable and unstable manifolds in $\Pi_\Dl$  ($\Dl \neq 0$).

Figure 2.6: The center-stable and center-unstable manifolds of $l_u$.

Figure 2.7: Persistent center-stable and center-unstable manifolds.

Figure 4.1: The graph of $\al = \al (\om)$.

\newpage
\bibliography{NLS}

\end{document}